\documentclass[a4paper,12pt]{article}

\usepackage{amsfonts} 
\usepackage{amsmath} 
\usepackage[american]{babel} 
\usepackage{amssymb} 
\usepackage{enumerate} 
\usepackage{theorem} 
\usepackage{graphicx} 
\usepackage{epsf} 
 
\newcommand{\eu}{{\bf E}^2} 
\newcommand{\eud}{{\bf E}^d} 
 
\newcommand{\fuso}{{\mathfrak H}}

\newcommand{\ds}{\displaystyle} 
\newcommand{\reach}{{\rm reach}} 
\newcommand{\unp}{{\rm Unp}} 
\newcommand{\tang}{{\rm Tan}} 
\newcommand{\nor}{{\rm Nor}} 
\newcommand{\norma}{\Vert}

\renewcommand{\sectionmark}[1]% 
   {\markright{\MakeUppercase{\thesection.\ #1}}} 
\newtheorem{teo}{Theorem}[section] 
\newtheorem{lm}[teo]{Lemma} 
\newtheorem{coro}[teo]{Corollary} 
\newtheorem{prop}[teo]{Proposition} 
\newtheorem{oss}[teo]{Remark} 
\newtheorem{Def}[teo]{Definition} 
 
\newenvironment{dimo}{{\it Proof.}}{$\square$\vspace{0.3cm}}

\begin{document} 
 
\title{On geometric properties of sets of positive reach in $\eud$} 
\author{Andrea Colesanti and Paolo Manselli} 
\date{} 

\maketitle 
\begin{abstract} 
\noindent
Some geometric facts concerning sets with positive
reach in $\eud$ are proved. For $x_1$ and $x_2$ in $\eud$ and $R>0$ let us denote by $\fuso(x_1,x_2,R)$ the 
intersection of all closed balls of radius $R$ containing $x_1$ and 
$x_2$. We prove that $\reach(K)\ge R$ if and only if for every $x_1,x_2\in K$ such that $\norma 
x_1-x_2\norma< 2R$, $\fuso(x_1,x_2,R)\cap K$ is connected. 
A corollary is that if $\reach(K)\ge R>0$ and $D$ is a closed ball of radius less than or equal to
$R$ (intersecting $K$) then $\reach(K\cap D)\ge R$. For $A\subset\eud$ and
$R>0$ we say that $A$ admits {\em $R$-hull} if there extsts $\hat A\supset A$, with
${\rm reach}(\hat A)\ge R$ and such that $\hat A$ is the minimal set (with respect to inclusion)
having these properties. A necessary and sufficient condition for a set
$A\subset\eud$ to admit a $R$-hull is provided.

\end{abstract}

\bigskip

 
\section{Introduction} 
Sets of positive reach were introduced by Federer in \cite{Federer}. This 
class of sets can be viewed as an extension of that of convex sets. It is well 
known that every point $x$ external to a closed convex set $C$ in $\eud$ admits a 
unique {\it projection} on $C$, i.e. a point which minimizes the distance from 
$x$ among all points in $C$. Sets of positive reach are those for which the 
projection is unique for the points of a parallel neighborhood of the set (and 
not necessarily for all external points). 
 
Along with their definition, Federer provided the main fundamental 
properties of sets of positive reach. Namely, the validity of global and 
local Steiner formulas and consequently the existence of curvature measures 
and many relevant properties of such measures. 
 
The study of properties of sets with positive reach has been continued by 
several authors and along various directions. Let us mention the contributions 
given by Z\"ahle \cite{Zhale} and Rataj and Z\"ahle \cite{Rataj-Zhale} on integral representation of curvature 
measures, the results by Hug \cite{Hug}, and Hug and the first author 
\cite{Colesanti-Hug} on singular points of sets with positive reach and the 
extensions of Steiner type formulas by Hug, Last and Weil 
\cite{Hug-Last-Weil}. Moreover, in \cite{Fu} Fu proved several interesting 
connections between sets of positive reach and semi-convex functions. 
 
As stated by Federer, closed convex sets represent a limit case of sets 
of positive reach, as the reach tends to $\infty$. The following question was 
at the origin of the research carried out in this paper. Is it possible to see 
(at least some of) the geometric properties of convex sets as limit case of 
suitable geometric properties of sets of positive reach? 
 
The first property that we analyse is the very definition of convex set: if 
$x_1$ and $x_2$ belong to a convex set $C$, then the segment joining them 
is entirely contained in $C$. In \S 3 we prove a 
possible counterpart of this fact for sets of positive reach. For two points 
$x_1$ and $x_2$ in $\eud$ and $R>0$ we denote by $\fuso(x_1,x_2,R)$ the 
intersection of all closed balls of radius $R$ containing $x_1$ and 
$x_2$. The set $\fuso(x_1,x_2,R)$ is a rugby ball-shaped set with cusps in $x_1$ and 
$x_2$; moreover for $R\to\infty$, $\fuso(x_1,x_2,R)$ tends to the segment 
with endpoints $x_1$ and $x_2$. Theorem \ref{T2.1} states that 
$\reach(K)\ge R$ if and only if for every $x_1,x_2\in K$ such that $\norma 
x_1-x_2\norma< 2R$, $\fuso(x_1,x_2,R)\cap K$ is connected. The proof of this 
result is geometric and does not require sophisticated techniques. 
As a corollary (see Theorem \ref{T2.2}) we have the following fact: if
$\reach(K)\ge R>0$ and $D$ is a closed ball of radius less than or equal to
$R$, intersecting $K$, then $\reach(K\cap D)\ge R$. The latter property can be seen as a
counterpart, for sets with positive reach, of the well-known fact that the
intersection of a convex set with an half-space is convex (if it is non-empty).

Next, we consider the following problem: given a set $A$ and a number $R>0$ 
is it possible to find the minimal set (with respect to inclusion) containing $A$ and having 
reach greater than or equal to $R$? The corresponding problem in the context of 
convexity ($R=\infty$) has an affirmative answer: every set admits a least convex cover, 
i.e. its convex hull. We will see through simple examples that this is not the 
case for arbitrary $A$ and $R$ and we will find necessary and sufficient conditions so that 
$A$ admits a minimal cover of reach greater than or equal to $R$. 
 
The paper is organized as follows: in \S 2 we introduce some notations; in \S 3 we 
prove Theorem \ref{T2.1} and some related results; in \S 4 we deal with the 
least cover with prescribed reach of a given set.

\section{Notations} 
Let $\eud$ be the $d$-dimensional Euclidean space; for $a,b\in\eud$, 
let $\Vert b-a\Vert$ be their distance and let $(\cdot,\cdot)$ denote 
the usual scalar product. 
 
If $A$ is a subset of $\eud$, then ${\rm int}(A)$, ${\rm cl}(A)$ and $A^c$ will 
denote the interior, the closure and the complement set of $A$, 
respectively. For $x_0\in\eud$ and $r>0$ we set 
$$ 
B(x_0,r)=\{x\in\eud\,:\,\Vert x-x_0\Vert<r\}\,,\quad{\rm and}\quad 
D(x_0,r)={\rm cl}(B(x_0,r))\,. 
$$ 
For $A\subset\eud$ and $a\in\eud$, the distance of $a$ from $A$ is given by 
$$ 
\delta_A(a)=\inf\{\Vert a-x\Vert\,:\, x\in A\}\,. 
$$ 
 
Let us recall the definition of {\em sets of positive reach}, introduced in 
 \cite{Federer}. Let $K\subset\eud$ be closed; let ${\rm Unp}(K)$ be the set 
 of points having a unique projection (or foot point) on $K$: 
$$ 
{\rm Unp}(K):=\{a\in\eud\,:\,\exists!\,x\in K\;{\rm s.t.}\;\delta_K(x)=\Vert a-x\Vert\}\,. 
$$ 
This definition implies the existence of a projection mapping 
$\xi_K\,:{\rm Unp}(K)\to K$ which assigns to $x\in{\rm Unp}(K)$ the unique 
point $\xi_K(x)\in K$ such that $\delta_K(x)=\Vert x-\xi_K(x)\Vert$. For a 
point $a\in K$ we set: 
$$ 
\reach(K,a)=\sup\{r>0\,:\, B(a,r)\subset\unp(K)\}\,. 
$$ 
The reach of $K$ is then defined by: 
$$ 
\reach(K)=\inf_{a\in K}\reach(K,a)\,, 
$$ 
and $K$ is said to be of positive reach if $\reach(K)>0$. 
 
If $K\subset\eud$ is compact and $x\in K$, the tangent and the normal spaces to 
$K$ at $a$ are: 
$$ 
\tang(K,a)=\{0\}\cup 
\left\{u\,:\,\forall\,\epsilon>0\;\exists\,b\in K\,{\rm s.t.}\, 
0<\Vert b-a\Vert<\epsilon\,,\, 
\left\Vert\frac{b-a}{\Vert b-a\Vert}-\frac{u}{\Vert u\Vert}\right\Vert<\epsilon\right\}\,, 
$$ 
$$ 
\nor(K,a)=\{v\,:\,(u,v)\le 0\,,\,\forall u\in\tang(K,a)\}\,. 
$$ 
Notice in particular that $\nor(K,a)$ is a closed convex cone. Let 
$\reach(K)>0$; for $a\in K$ we set: 
$$ 
P_a=\{v\,:\,\xi_K(a+v)=a\}\,,\quad 
Q_a=\{v\,:\,\delta_K(a+v)=\Vert v\Vert\}\,. 
$$ 
 
\section{Characterization and geometrical properties of sets with positive reach } 
 
The following definition will be useful later. 
 
\begin{Def}\label{D2.1}  Let $a,b\in\eud$, $a\ne b$, $R>0$, $\Vert a-b\Vert<2R$. Let 
$$ 
{\cal D}(a,b,R)=\{D(x,R)\,:\,\norma a-x\norma,\norma b-x\norma\le R\}\,. 
$$ 
We set 
$$ 
\fuso(a,b,R)=\bigcap_{D\in{\cal D}(a,b,R)} D\,. 
$$ 
\end{Def} 
 
It is clear from the definition that $\fuso(a,b,R)$ is a compact convex set, 
containing $a$ and $b$. The boundary of $\fuso(a,b,R)$ is 
obtained rotating an arc of circle of radius $R$ joining $a$ and $b$, about the 
line through $a$ and $b$. 
 
\begin{lm}\label{L2.3} Let $a,b\in\eud$ be such that $0<\Vert b-a\Vert< 2R$ 
  where $R>0$. If $c,d\in\fuso(a,b,R)$, then $\fuso(c,d,R)\subset\fuso(a,b,R)$. 
\end{lm} 
\noindent
\begin{dimo} If $D\in{\cal D}(a,b,R)$, then $c,d\in D$ so that $D\in{\cal D}(c,d,R)$. 
The conclusion follows from Definition \ref{D2.1}. 
\end{dimo} 
 
A set is convex if and only if given any two points belonging to it, it contains 
the line segment joining them. 
In this section we prove (see Theorem \ref{T2.1}) 
a characterization of sets of positive reach that somehow 
resembles the above characterization of convex sets. 
 The proof of this result requires various lemmas. 
%\begin{lm}\label{L2.1} 
%Let $K\subset\eud$ be compact, $\reach(K)\ge R>0$. If 
%  $a\in\partial K$ then 
%$$ 
%P_a\cap B(0,R)=Q_a\cap B(0,R)=\nor(K,a)\cap B(0,R)\,. 
%$$ 
%\end{lm} 
%\begin{dimo} By Theorem 4.8 (2) in \cite{Federer} we have 
%\begin{equation} 
%\label{E2.1} 
%P_a\subset Q_a\subset\nor(K,a)\,. 
%\end{equation} 
%Let $w\in\nor(K,a)\cap B(0,R)$, $\Vert w\Vert=r<R$. By Theorem 4.8 (12) in 
%\cite{Federer} there exists $v\in P_a$ such that $\Vert v\Vert=r$ and 
%$w=\lambda v$ for some $\lambda\ge 0$. This implies $\lambda =1$ and $w=v$ 
%so that $w\in P_a$. This shows that $\nor(K,a)\cap B(0,R)\subset P_a\cap 
%B(0,R)$, and, by (\ref{E2.1}), the proof is concluded. 
%\end{dimo} 
The next proposition is Theorem 4.8 (7) of \cite{Federer}. 
 
\begin{prop}\label{L2.2} Let $K\subset\eud$ be closed, $x\in\unp(K)$ and 
  $\reach(K,\xi_K(x))>0$. Then, for every $b\in K$ 
\begin{equation} 
\label{E2.2} 
(x-\xi_K(x),\xi_K(x)-b)\ge- 
\frac{\Vert\xi_k(x)-b\Vert^2\,\Vert x-\xi_K(x)\Vert}{2\,\reach(K,\xi_K(x))}\,. 
\end{equation} 
\end{prop} 
 
Let $R>0$ and $a,b\in\eud$ be such that 
$0<\norma a-b\norma<2R$. We define the cone 
$$ 
{\cal C}(a,b,R)=\left\{v\ne0 \,:\,\left(\frac{v}{\norma v\norma},\frac{b-a}{\norma 
      b-a\norma}\right)>\frac{\norma b-a\norma}{2R}\right\}\,. 
$$ 
A geometric version of the above proposition follows. 
 
\begin{coro}\label{C2.1} Let $K$ be a closed subset of $\eud$ such that 
  $\reach(K)\ge R>0$. Let $x\in\unp(K)\setminus K$, $a=\xi_K(x)\in\partial K$ and $b\in K$ 
such that $0<\norma a-b\norma<2R$. Then 
$$ 
x-a\notin{\cal C}(a,b,R)\,. 
$$ 
\end{coro} 
 
We proceed with some geometric considerations in the plane. 
Given $v$ and $w$ vectors in $\eu$, $v,w\ne0$, we set 
$$ 
{\cal S}(v,w)=\{z\,:\,z=tv+\tau w\,,\,t,\tau> 0\}\,. 
$$ 
 
\begin{oss}\label{O2.1} Let $R>0$ and $z_1,z_2,z_3,z_4\in\eu$ be such that 
$$ 
\norma z_1-z_2\norma=\norma z_2-z_3\norma= 
\norma z_3-z_4\norma=\norma z_4-z_1\norma=R\,,\quad 0<\norma z_1-z_3\norma<2R\,. 
$$ 
We have 
$$ 
{\cal C}(z_1,z_3,R)={\cal S}(z_2-z_1,z_4-z_1)\,. 
$$ 
\end{oss} 
 
\begin{lm}\label{L2.5} Let $R>0$, $b_1,b_2\in\eu$ with $0<\Vert b_1-b_2\Vert<2R$, 
  $\Gamma_j=\partial B(b_j,R)$, $j=1,2$, $b_3,b_4\in\eu$ such that 
  $\{b_3,b_4\}=\Gamma_1\cap\Gamma_2$. Let 
\begin{enumerate} 
\item[(i)] $\Sigma\subset\Gamma_1$ be the closed 
  arc joining $b_3$ and $b_4$ of smaller length; 
\item[(ii)] $\Sigma'\subset\Gamma_1$ 
  be the closed arc having length $\pi R$ and such that 
  $\Sigma\cap\Sigma'=\{b_4\}$. 
\end{enumerate} 
For every $a\in B(b_4,R)\setminus 
  D(b_3,R)$ there exist $c\in\Sigma$, $c\ne b_3$, $c\ne b_4$, and 
  $c'\in\Sigma'$, uniquely determined, such that 
$$ 
\norma b_1-c' \norma=\norma c'-a   \norma=  
\norma a-c    \norma=\norma c-b_1  \norma=R\,. 
$$ 
\end{lm} 

\vskip0.5cm
\begin{center}
\includegraphics[scale=.4]{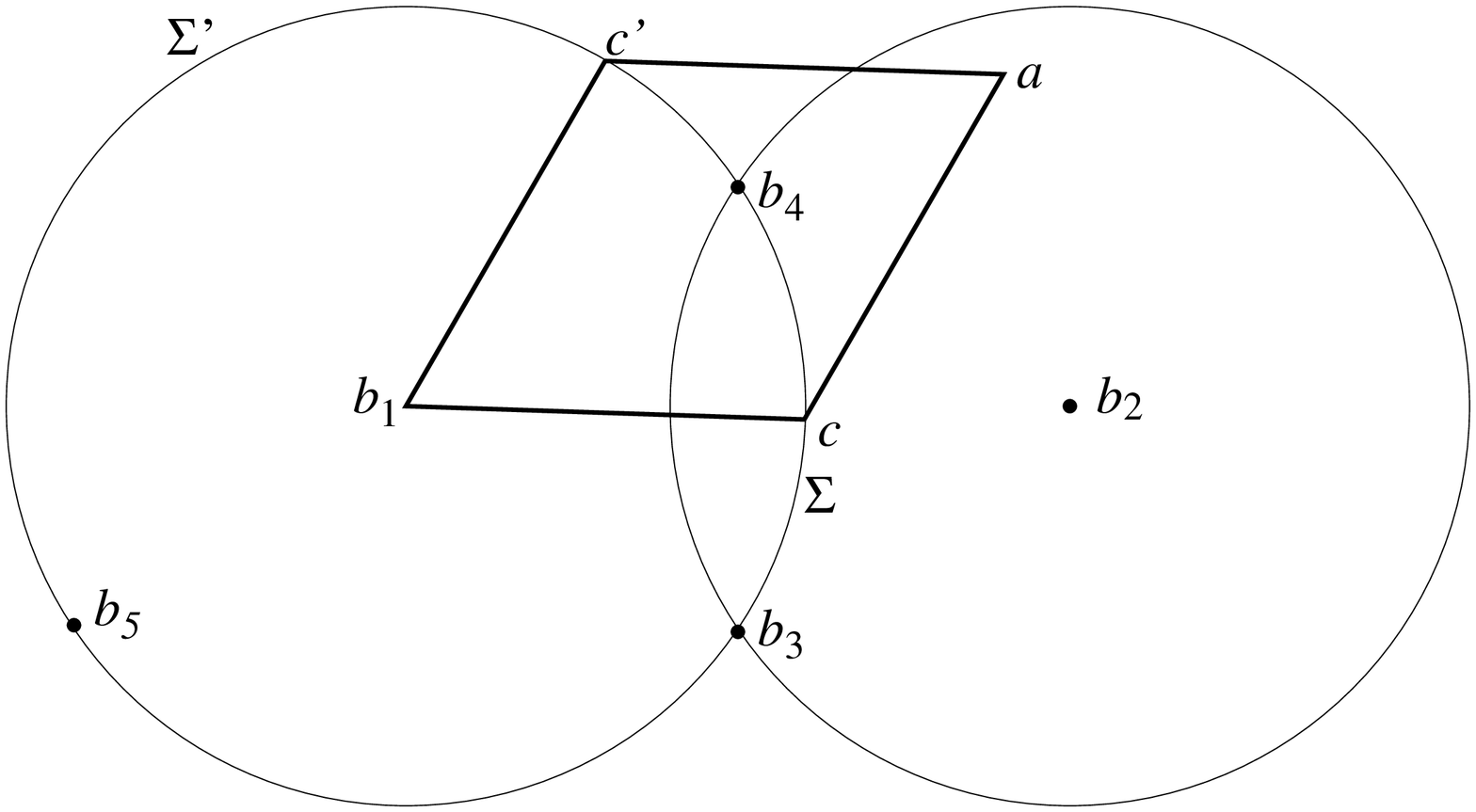} 
\vskip0.5cm
{\sc Figure 1}
\end{center}

\noindent 
\begin{dimo} We have $\Vert a-b_3\Vert>R$ and $\Vert a-b_4\Vert< R$. 
  Let us notice that $b_3$ and $b_4$ are the endpoints of $\Sigma$. By continuity, there 
  exists $c\in\Sigma$ such that $\Vert a-c\Vert=R$. Let $b_5$ be the endpoint 
  of $\Sigma'$ which does not coincide with $b_4$; we have $\Vert 
  b_4-b_5\norma=2R$ and $\Vert a-b_5\Vert+\norma a-b_4\norma\ge\norma  
  b_4-b_5\norma=2R$; thus $\norma a-b_5\norma\ge2R-\norma a-b_4\norma> R$. By 
  continuity, there exists $c'\in\Sigma'$ such that $\norma a-c'\norma=R$. The  
  points $c$ and $c'$ are uniquely determined as intersection of 
  $\Gamma_1$ and $\partial B(a,R)$. 
\end{dimo}

\begin{lm}\label{L2.6} Let $R>0$, $b_1,b_2\in{\bf E}^2$, $0<\norma b_1-b_2\norma<2R$, $B_i=B(b_i,R)$, 
  $\Gamma_i=\partial B_i$, $i=1,2$. Let $b_3,b_4$ be such that 
  $\{b_3,b_4\}=\Gamma_1\cap\Gamma_2$, $B_i=B(b_i,R)$, 
  $i=3,4$. Assume that $a\in B_3\cup B_4\setminus\fuso(b_1,b_2,R)$ and $c_i,c'_i$ are 
  such that 
$$ 
\norma b_i-c'_i \norma=\norma c'_i-a   \norma= 
\norma a-c_i    \norma=\norma c_i-b_i  \norma=R\,, 
\quad{\rm for}\;i=1,2\,, 
$$ 
and let $S_i=S(c_i-a,c'_i-a)$, for $i=1,2$. Then: 
\begin{equation} 
\label{E2.4} 
S_1\cup S_2\supset S(b_2-a,b_1-a)\,. 
\end{equation} 
In particular 
\begin{equation}\label{E2.4b} 
\frac{1}{2}(b_1+b_2)\in{\rm int}(S_1\cup S_2)\,. 
\end{equation} 
\end{lm} 

\vskip0.5cm
\begin{center}
\includegraphics[scale=.4]{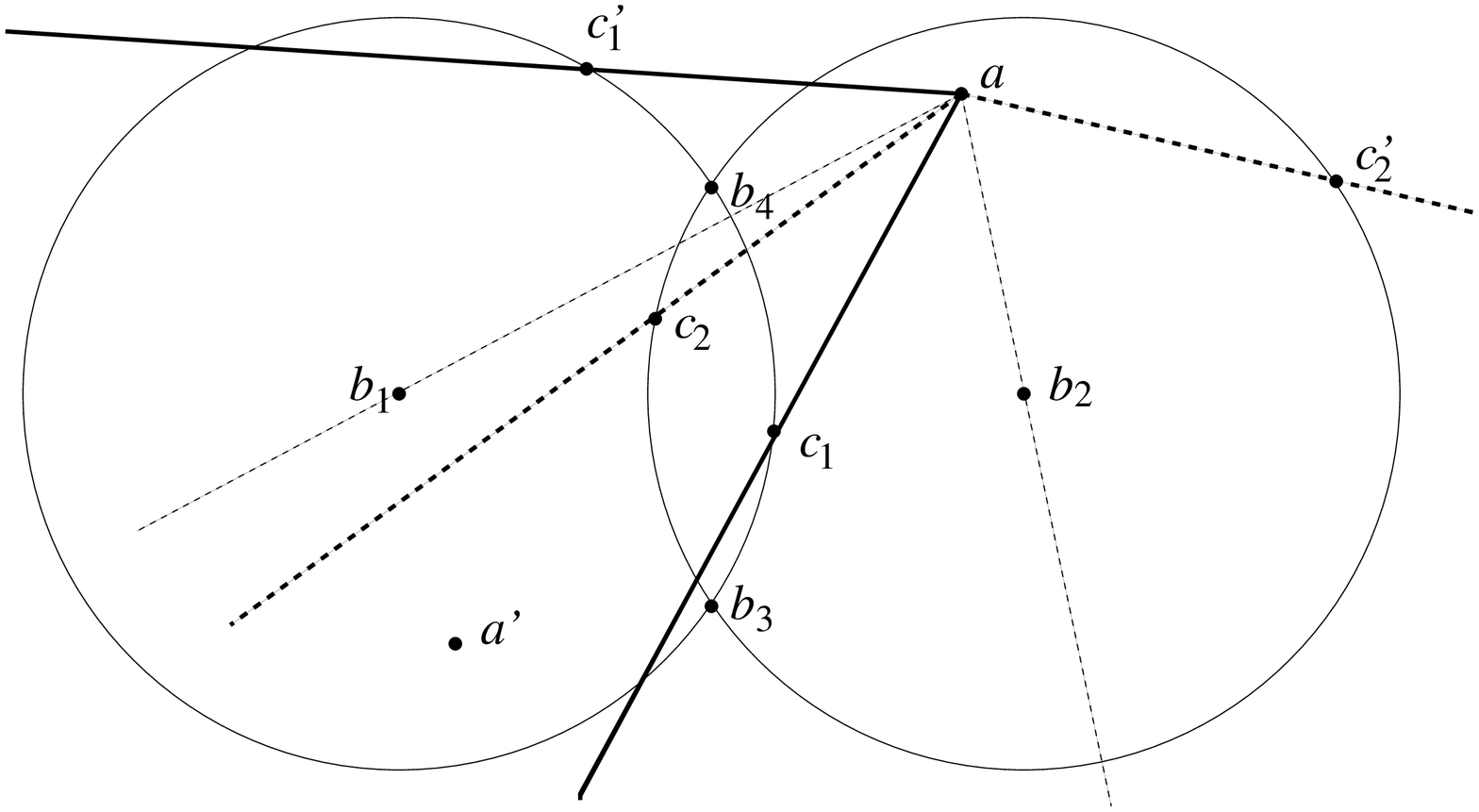} 
\vskip0.5cm
{\sc Figure 2}
\end{center}

\noindent 
\begin{dimo} 
$S_1, S_2$ and $S(b_2-a,b_1-a)$ are open convex sectors with apex in $a$; 
moreover $b_i-a\in S_i$ for $i=1,2$ so that $\{b_1-a,b_2-a\}\subset S_1\cup 
S_2$. 
Let $\Sigma_1=\Gamma_1\cap D(b_2,R)$ and $\Sigma_2=\Gamma_2\cap D(b_1,R)$. By 
Lemma \ref{L2.5} we may assume that $c_i\in\Sigma_i\setminus\{b_3,b_4\}$ for $i=1,2$. This in turn 
implies $\Vert c_1-c_2\Vert<2R$ (as $c_1,c_2\in\fuso(b_3,b_4,R)$). Hence it is uniquely determined 
$a'\ne a$ such that $\{a,a'\}=\partial B(c_1,R)\cap\partial B(c_2,R)$. 
The straight line through $a$ and $a'$ bounds two open half-planes such that $b_2$ and $c_1$ 
(resp. $b_1$ and $c_2$) are in the same half-plane. Thus 
\begin{equation} 
\label{E2.5} 
a'-a\in S_1\cap S_2\ne\emptyset\,. 
\end{equation} 
This implies that $S_1\cup S_2$ is a convex cone and,
since it contains $b_1$ and $b_2$, (\ref{E2.4}) follows.
\end{dimo}

\begin{teo}\label{T2.1} If $K\subset\eud$ is closed then $\reach(K)\ge R>0$ if 
  and only if for every $b_1,b_2\in K$, $\Vert b_1-b_2\Vert<2R$, $K\cap 
  \fuso(b_1,b_2,R)$ is connected. 
\end{teo} 

\noindent 
\begin{dimo} Let us assume that $\reach(K)\ge R>0$.
  By contradiction, assume that 
  $K':=K\cap\fuso(b_1,b_2,R)$ 
  is not connected; then there exist $K_1,K_2\subset K'$, 
  closed, such that $K'=K_1\cup K_2$ and $K_1\cap 
  K_2=\emptyset$. By compactness, there exist $c_i\in K_i$ for $i=1,2$ 
  such that 
$$ 
\rho:=\Vert c_1-c_2\Vert=\inf\{\Vert x-y\Vert\,:\,x\in K_1\,,\, y\in 
K_2\}>0\,. 
$$ 
As $c_1,c_2\in\fuso(b_1,b_2,R)$, $\rho\le R$. We have 
$$ 
B(c_1,\rho)\cap B(c_2,\rho)\cap K'=\emptyset\,. 
$$ 
On the other hand it is easy to check that 
$\fuso(c_1,c_2,R)\subset[B(c_1,\rho)\cap B(c_2,\rho)]\cup\{c_1,c_2\}$. By 
Lemma \ref{L2.3}, $\fuso(c_1,c_2,R)\subset\fuso(b_1,b_2,R)$, so that 
\begin{equation} 
\label{E2.6} 
\fuso(c_1,c_2,R)\cap K'=\{c_1,c_2\}\,. 
\end{equation} 
In particular, $c:=\ds{\frac{c_1+c_2}{2}}\notin K$; as 
$\delta_K(c)<R$, $c\in\unp(K)\setminus K$. Let $c_3=\xi_K(c)\in\partial 
K$. Notice that if $c_3\in\fuso(c_1,c_2,R)$ then either $c_3=c_1$ or $c_3=c_2$ so 
that $\delta_K(c)=\Vert c-c_1\Vert=\Vert c-c_2\Vert$ in contradiction with 
$c\in\unp(K)$. Consequently, $c_3\in K\setminus\fuso(c_1,c_2,R)$. 
We also observe that, for $i=1,2$, 
$$ 
\norma c_i-c_3\norma\le\norma c_i-c\norma+\norma c-c_3\norma<2R 
$$ 
as $\norma c-c_3\norma=\delta_K(c)<R$. 
We recall the definitions of the cones: 
$$ 
{\cal C}_i(c_3,c_i,R)=\left\{v\ne0 \,:\,\left(\frac{v}{\norma v\norma},\frac{c_i-c_3}{\norma 
      c_i-c_3\norma}\right)>\frac{\norma c_i-c_3\norma}{2R} 
\right\}\,,\quad i=1,2\,. 
$$ 
By Corollary \ref{C2.1} we have that 
\begin{equation} 
\label{E2.7} 
c-c_3\notin{\cal C}_1\cup{\cal C}_2\,. 
\end{equation} 
Apply Remark \ref{O2.1} and Lemma \ref{L2.6} to the (uniquely determined) 2-dimensional plane 
containing $c,c_1,c_2,c_3$ to obtain a contradiction with (\ref{E2.7}). 
 
Vice-versa, assume that for every $b_1,b_2\in K$, $\Vert b_1-b_2\Vert<2R$, the set 
  $K\cap 
  \fuso(b_1,b_2,R)$ is connected. If, by contradiction, $\reach(K)< R$, 
  then there exists $x\in K^c$ such that $\delta_K(x)=r<R$ and 
  $\norma x-b_1\norma=\norma x-b_2\norma=r$ for some $b_1,b_2\in K$, 
  $b_1\ne b_2$. As $\norma b_1-b_2\norma<2R$, $\fuso(b_1,b_2,R)\cap K$ is 
  connected. On the other hand, $r<R$ implies that $\fuso(b_1,b_2,R)\subset 
  B(x,r)\cup\{b_1,b_2\}$ so that there exists $b\in K\cap B(x,r)$ i.e. a 
  contradiction. 
\end{dimo} 
 
\begin{oss} If $\reach(K)\ge R$ and $b_1,b_2\in K$ are such that 
$\norma b_1-b_2\norma=2R$, then $K\cap\fuso(b_1,b_2,R)$ is not necessarily 
connected. Any set consisting of two points at distance $2R$ is an example. 
\end{oss}

\begin{teo}\label{T2.2} Let $K$ be a closed set such that $\reach(K)\ge R>0$. If $D$ is a 
closed ball of radius less than or equal to $R$ then $\reach(K\cap 
D)\ge R$. 
\end{teo} 
 
\noindent 
\begin{dimo} The argument is similar to the one used in the second part of
  the proof of Theorem \ref{T2.1}. Let $a\in(K\cap D)^c$ such that $r=\delta_{K\cap D}(a)<R$; let us show 
that $a\in\unp(K\cap D)$. Assume by contradiction that there exist 
$b_1,b_2\in(K\cap D)$ such that $b_1\ne b_2$ and 
$\norma a-b_1\norma=\norma a-b_2\norma=r$. In particular $\norma b_1-b_2\norma<2R$. 
Clearly, $\fuso(b_1,b_2,R)\subset D$; consequently, by Theorem \ref{T2.1}, 
$\fuso(b_1,b_2,R)\cap(K\cap D)$ is connected. Also, notice that 
$$ 
(\fuso(b_1,b_2,R)\setminus\{b_1,b_2\})\subset B(a,r)\,. 
$$ 
Then there exists $b'\in K\cap D$ such that $\norma a-b'\norma<r$, i.e. a 
contradiction. 
\end{dimo} 
 
\begin{coro} If $\reach (K)\ge R>0$, $a,b\in K$, $\norma a-b\norma\le2R$, then 
$\reach(K\cap\fuso(a,b,R))\ge R$. 
\end{coro} 
 
It is well known that, if $ K $ is a closed convex set in $\eud$ and 
$ H $ is an open half space, satisfying $ H \cap K = \emptyset, $ 
then $ \partial H \cap K $ is either empty or a convex subset of $ 
\partial H. $ Let us show that a similar property holds for sets of 
reach $ \ge R>0.$ 
 
\begin{Def}\label{D2.2} 
Let  $S$  be a sphere of radius $  R>0 $ in $\eud$;  let $ K $ be 
a closed subset of $S$. We say that $ K $ is convex in $ S $ 
if $ x_1 \in K,\; x_2 \in K, $ dist($ x_1, x_2$) $< 2R $ imply 
that the arc of great circle of $S$ joining  $ x_1  $ and  $ x_2, $ and having
smaller length, is contained in $K$. 
\end{Def} 
 
\begin{teo}\label{T2.3} Let $K$ be a closed  set in $\eud$ 
and $\reach (K)\ge R>0$. Let $ B $ be an open ball of radius $R$ 
satisfying $ B \cap K = \emptyset$. Then $ \partial B \cap K $ is 
either empty or a convex subset of $ \partial B $. 
\end{teo} 
 
\noindent
\begin{dimo} 
Theorem \ref{T2.2} implies that $ (B\cup \partial B) \cap K $ $ = 
 \partial B \cap K $ has reach $ \ge R. $ Then, by theorem 
 \ref{T2.1}, if $b_1,b_2 \in K\cap \partial B$, $\Vert b_1-b_2\Vert<2R$, then $K\cap 
 \partial B \cap \fuso(b_1,b_2,R)$ is connected. Now 
$K\cap  \partial B \cap \fuso(b_1,b_2,R)$ is exactly the arc of great circle
 of $\partial B$, joining  $ b_1$ and  $ b_2$ and having smaller length. 
\end{dimo}

\section{On the $R$-hull of a set} 
 
Let $A$ be a subset of $\eud$ and let $R>0$. In this section we analyze 
the problem of finding $K$ such that $\reach(K)\ge R$, $K\supset A$ and $K$ is 
the {\em minimal} set (with respect to inclusion) having these properties. In 
other words we look for a sort of hull of reach $R$ of $A$. Intuitively, when 
$R=\infty$ we are dealing with the convex hull of $A$ which exists for 
every $A$. On the other hand, for finite $R>0$ not every set $A$ 
admits a hull of reach $R$ (see the examples below). Our aim is to give necessary and 
sufficient conditions for $A$ to have this property (see Theorems \ref{T4.1}
and \ref{T4.2}).  
 
\begin{Def}\label{D4.1} 
Let $A\subset\eud$, $R>0$. We say that $A$ admits a $R$-hull if 
  there exists $\hat A\subset\eud$ such that: 
\begin{enumerate} 
\item[(i)] $A\subset\hat A$; 
\item[(ii)] $\reach(\hat A)\ge R$; 
\item[(iii)] if $\reach(K)\ge R$ and $A\subset K$, then 
  $\hat A\subset K$. 
\end{enumerate} 
If such a set exists, we call it the $R$-hull of $A$. 
\end{Def} 
 
\noindent 
{\bf Example 1.} For an arbitrary $R>0$ we may construct 
an example of set which does not admit a $R$-hull. 
Let $n=2$  and $A=\{a,b\}$ with $\norma a-b\norma=R/2$. Assume by contradiction 
that there exists the $R$-hull of $A$, and denote it by $\hat A$. Let $\hat A_1$ be the 
closed line segment joining $a$ and $b$: $\reach(\hat A_1)=\infty$ so that 
$\hat A_1\supset \hat A$. Let $\Gamma$ be a circle of radius $R$ passing through $a$ and 
$b$ and let $\hat A_2\subset\Gamma$ be the closed arc of smaller length joining $a$ 
and $b$. We have $\reach(\hat A_2)=R$ so that $\hat A_2\supset \hat A$. As 
$\hat A_1\cap \hat A_2=A$, 
we must have $\hat A=A$; on the other hand $\reach(A)=R/2$ so we have a 
contradiction. 
 
\bigskip 
 
\noindent 
{\bf Example 2.} In $\eud$ consider a half-line $L$ with end-point in the origin. For 
every $i=1,2,\dots$, let $a_i$ be the point of $L$ such that 
$\norma a_i\norma=1/i$. The set $A=\{a_1,a_2,\dots\}$ does not admit a $R$-hull for any 
$R\in(0,\infty)$. 
 
\bigskip 
 
For an arbitrary set $A\subset\eud$ and $R>0$, we set 
$$ 
A'_R=\{x\in\eud\,:\,\delta_A(x)\ge R\}\,. 
$$ 
The proof of the following proposition is an easy application of Theorem \ref{T2.1}.  
\begin{prop} Let $A\subset\eud$, $R>0$; $\reach(A'_R)\ge R$ 
    if and only if for every $a$ and $b$ such that 
     $\delta_A(a),\delta_A(b)\ge R$ and $B(a,R)\cap B(b,R)\ne\emptyset$, 
    there exists a continuous arc $\Gamma$ joining 
    $a$ and $b$, $\Gamma\subset\fuso(a,b,R)$, 
such that $\delta_A(x)\ge R$ for every $x\in\Gamma$. 
\end{prop} 
 
\begin{lm}\label{L4.1} 
Let $K\subset\eud$, then 
\begin{enumerate} 
\item[(i)] $K\subset (K'_R)'_R\subset\{z\in\eud\,:\,\delta_K(z)<R\}$, 
\item[(ii)]  if $\reach(K)\ge R>0$ then $\reach(K'_R)\ge R$ and $K=(K'_R)'_R$. 
\end{enumerate} 
\end{lm} 
\noindent
\begin{dimo} If $x\in K$, then $\norma x-y\norma\ge R$ for every $y\in K'_R$ so that 
$\delta_{K'_R}(x)\ge R$ and $x\in(K'_R)'_R$. On the other hand, if $z\in(K'_R)'_R$ then 
$z\notin K'_R$ so that $\delta_K(z)<R$. Claim (i) is proved. 
 
For $s\ge 0$ set $K'_s=\{x\in\eud\,:\,\delta_K(x)\ge s\}$. 
Corollary 4.9 in \cite{Federer} implies that $\reach(K'_{R-1/i})\ge R-1/i$ 
for every $i=1,2,\dots$. Moreover, the sequence $K'_{R-1/i}$ converges to 
$K'_R$ in the Hausdorff metric. On the other hand, the by Remark 4.14 in 
\cite{Federer}, for every $\epsilon>0$ the family 
$$ 
\{A\subset\eud\,:\,\reach(A)\ge R-\epsilon\} 
$$ 
is closed with respect to the Hausdorff metric. Then $\reach(K'_R)\ge R-\epsilon$ 
for every $\epsilon>0$. Now let us prove that if $\reach(K)\ge R$ then 
$(K'_R)'_R\setminus K$ is empty. Let $z\in(K'_R)'_R\setminus K$; (i) implies that 
$z\in\unp(K)$. Let $x=\xi_K(z)$ and $y_t=x+t\frac{z-x}{\norma z-x\norma}$,
$t\ge 0$. Note that $\frac{z-x}{\norma z-x\norma}\in{\rm Nor}(K,x)$ so that,
by claim (12) of Theorem 4.8 of \cite{Federer}, if $0<t<R$, then $\delta_K(y_t)=t$ and by continuity
$\delta_K(y_R)=R$. Then $y_R\in K'_R$ and  
$\norma z-y_R\norma<R$, i.e. a contradiction. 
\end{dimo}

\begin{teo}\label{T4.1} Let $A\subset\eud$ and $R>0$. If $\reach(A'_R)\ge R$ then 
$A$ admits $R$-hull $\hat A$ and 
$$ 
\hat A=(A'_R)'_R\,. 
$$ 
\end{teo} 
\noindent 
\begin{dimo} Let $A_1=(A'_R)'_R$; we prove 
that $A_1$ is the $R$-hull of $A$. The inclusion $A\subset A_1$ is part (i) in Lemma 
\ref{L4.1}. 
By the same lemma, as $\reach(A'_R)\ge R$ we have $\reach(A_1)\ge R$. It remains to 
show that $A_1$ satisfies (iii) in Definition \ref{D4.1}. Let $K$ be such 
that $K\supset A$ and $\reach(K)\ge R$. Then $K'_R\subset A'_R$ and, by Lemma 
\ref{L4.1}, $K=(K'_r)'_R\supset(A'_R)'_R=A_1$. 
\end{dimo} 
 
\begin{coro}Let $A\subset\eud$ and $R>0$. If for every $a$ and $b$ such that 
$\delta_A(a),\delta_A(b)\ge R$ and $\norma a-b\norma<2R$, there exists a continuous 
arc $\Gamma$,  joining $a$ and $b$ such that $\delta_A(x)\ge R$ for every $x\in\Gamma$, 
$\Gamma\subset\fuso(a,b,R)$, then $A$ admits $R$-hull $\hat A$ and 
$$ 
\hat A=(A'_R)'_R\,. 
$$ 
\end{coro} 
 
\begin{teo}\label{T4.2} Let $K\subset\eud$ and $R>0$. Assume that $K$ admits $R$-hull 
$\hat K$. Then $\reach(K'_R)\ge R$. 
\end{teo} 
\noindent 
\begin{dimo} We argue by contradiction. By using Theorem \ref{T2.1}, 
there exist $b_1$ and $b_2\in K'_R$ satisfying $\norma b_1-b_2\norma<2R$ and such that 
$\fuso(b_1,b_2,R)\cap K'_R$ is not connected. Then, as we saw in the proof of 
Theorem \ref{T2.1}, there exist $c_1$ and $c_2\in K'_R$ such that 
\begin{equation} 
\label{4.1} 
\fuso(c_1,c_2,R)\cap K'_R=\{c_1,c_2\}\,. 
\end{equation} 
For $j=1,2$ we have $\reach(B(c_j,R)^c)=R$ and $B(c_j,R)^c\supset K$ thus 
$B(c_j,R)^c\supset\hat K$. This implies in particular that $c_1$, $c_2\in(\hat 
K)'_R$. As $\reach(\hat K)\ge R$, by Lemma \ref{L4.1}, 
$\reach(\hat K'_R)\ge R$, then $\fuso(c_1,c_2,R)\cap\hat K'_R$ is connected. 
Let $a\in[\fuso(c_1,c_2,R)\setminus\{c_1,c_2\}]\cap\hat K'_R$. We have 
$B_a(R)\cap K\subset B_a(R)\cap\hat K=\emptyset$ then $a\in K'_R$ which contradicts 
(\ref{4.1}). 
\end{dimo} 
 
From the above theorem another connection between convex sets and sets of
positive reach can be deduced. The convex hull of a closed set $C$ is the intersection of all the closed half-spaces 
containing $C$. Let us prove that if $K$ admits $R$-hull $\hat K$, then $\hat K$ is the 
intersection of the complement sets of all open balls that do not meet $K$. 
Note that for an arbitrary, non-empty, subset $K$ of $\eud$ we have 
$$ 
(K'_R)'_R=\bigcap_{\delta_K(x)\ge R}B_x(R)^c\,. 
$$ 
This remark and Theorem \ref{T4.2} lead to the following result. 
 
\begin{coro} Let $K\subset\eud$, $R\ge 0$. Assume that $K$ admits an $R$-hull 
$\hat K$. Then 
$$ 
\hat K=\bigcap_{\delta_K(x)\ge R}B_x(R)^c\,. 
$$ 
\end{coro}

\bigskip

{\sc Andrea Colesanti}, Dipartimento di Matematica 'U. Dini', Viale Morgagni 67/a, 50134 
Firenze, Italy. E-mail: {\tt colesant@math.unifi.it}

{\sc Paolo Manselli}, Dipartimento di Matematica e Applicazioni per l'Architettura, via dell'Agnolo 14, 50122 Firenze, 
Italy. E-mail: {\tt manselli@unifi.it}

\end{document}